\documentclass[leqno]{macrorend}
\usepackage{mathrsfs}
\volumeyear{xxxx}\yearnumber{x}\volumenumber{xx}

\begin{document}

% select a language among: english, french, italian
%
\selectlanguage{english}

% if you dont have footnote, cancel \footnotemark[1]
% separate several authors with ' - '
%
\articolo[Periodic orbits in Hill lunar problem]
{Doubly-symmetric periodic orbits in the spatial Hill's lunar problem
with oblate secondary primary}{Xingbo Xu\footnotemark[1]}

% text of the note (number from 1 to 9)
%
\footnotetext[1]{xuxingbo25@hotmail.com, the author is supported 
by the National Nature Science Foundation of China(NSFC)
with the Grant No. 11703006.}

\begin{abstract}
  In this article we consider the existence of a family of doubly-symmetric 
  periodic orbits in the spatial circular Hill's lunar problem, in which
  the secondary primary at the origin is oblate. The existence is 
  shown by applying a fixed point theorem to the equations with periodical conditions
  expressed in Poincar\'{e}-Delaunay elements for the double symmetries
  after eliminating the short periodic effects in the first-order perturbations
  of the approximated system.

\end{abstract}

\section{Introduction}

   The classical spatial Hill's lunar problem is a limiting case derived from the
   spatial \emph{circular restricted three-body problem} (CRTBP) \cite{Meyer1982}.
   Consider a modified version of the classical CRTBP, in such a problem,
   one small primary $\texttt{M}_{2}$ is an oblate Maclaurin ellipsoid and the other primary $\texttt{M}_{1}$
   is a standard spheroid, so $\texttt{M}_{1}$ can be considered as 
   a mass point while $\texttt{M}_{2}$ has a shape.
   Set the masses of $\texttt{M}_{1}$ and $\texttt{M}_{2}$ to be 
   $m_{1}$ and $m_{2}$, respectively. 
   Denote the radius of the equator of $\texttt{M}_{2}$ as $a_{e}$, the polar radius $b_{e}$.
   The relative position of $\texttt{M}_{1}$ relative to $\texttt{M}_{2}$ is $\vec{r}$, and denote $r$ 
    be the length of the vector $\vec{r}$.    
   In order to make the differential system about $\ddot{\vec{r}}$ integrable, 
   suppose $\texttt{M}_{1}$ moves on the equator plane of $\texttt{M}_{2}$. 
   The angular velocity $\omega_{e}$ of the relative circular motion of the two primaries
   can be calculated by
   \begin{align}
     \omega_{e}^{2}=\texttt{G}(m_{1}+m_{2})
     \left(1-\sum_{n=1}^{\infty}(2n+1)\frac{a_{e}^{2n}}{r^{2n+2}}J_{2n} P_{2n}(0)\right).
   \end{align}
   where $\texttt{G}$ is the universal gravitational constant, $\{J_{2n}\}$ are the even zonal harmonic 
   coefficients, and $P_{n}(x)$ is the n-th Legendre polynomial. 
   By the symbol computations, $J_{2n}$  for the Maclaurin ellipsoid can be calculated as 
   \begin{align}
       J_{2n}=(-1)^{n-1}\left(1-\frac{b_{e}^{2}}{a_{e}^{2}}\right)^{n}
       \prod_{k=1}^{n}\frac{2k-1}{2k+3}=\frac{2n-1}{2n+3}\left(1-\frac{b_{e}^{2}}{a_{e}^{2}}\right)J_{2n-2} ,
        \quad n\geq 1.       
   \end{align}
   Set the total masses as the mass unit, the radius and the period of relative circular motion 
   of the two primaries as the distance unit, and $2\pi$, respectively. In such units, the universal
   gravitational constant $\texttt{G}=1$, and the masses of $\texttt{M}_{1}$ and $\texttt{M}_{2}$ are 
   $1-\mu$ and $\mu$ respectively. Set the origin at the center of $\texttt{M}_{2}$, the positive $x$
   axis is the direction from $\texttt{M}_{2}$ to $\texttt{M}_{1}$, the position
   $\vec{r}=\mathbf{x}=(x_{1},x_{2},x_{3})^{T}$ is a column vector. 
   In the uniform rotating frame, the conjugate momentum of $\mathbf{x}=\vec{r}$ is 
   $\mathbf{y}=(\dot{x}_{1}-x_{2},\dot{x}_{2}+x_{1},\dot{x}_{3})^{T}$. The Hamiltonian for the motion
   of the infinitesimal body moving nearby $\texttt{M}_{2}$ can be written as
   \begin{align}
     H_{a}  &=\frac{1}{2}\left[y_{1}^{2}+(y_{2}-1+\mu)^{2}+y_{3}^{2}\right]
     -[(x_{1}-1+\mu)(y_{2}-1+\mu)-x_{2}y_{1}]
     -U_{1}-U_{2} ,
      \nonumber   \\
     U_{1} &=\frac{1-\mu}{\sqrt{1+r^{2}-2x_{1}}}
     =(1-\mu)\sum_{n=0}^{\infty}P_{n}(\frac{x_{1}}{r})r^{n} ,  \nonumber  \\
     U_{2} &=\frac{\mu}{r}
     -\frac{\mu}{r}\sum_{n=1}^{\infty}\left(\frac{a_{e}}{r}\right)^{2n}J_{2n}P_{2n}(\frac{x_{3}}{r}) .
   \end{align}
   Expand the Hamiltonian, neglect the constant term.
    Make symplectic scaling, $\mathbf{x}\rightarrow \mu^{1/3}\mathbf{x}$,
    $\mathbf{y}\rightarrow \mu^{1/3}\mathbf{y}$, at the mean while, set
    \begin{align*}
    J_{2n}=-C_{2n}, \quad 
     r \rightarrow  \mu^{1/3}r ,  \quad
     a_{e}\rightarrow  \mu^{1/3}a_{e}, \quad 
     b_{e}\rightarrow  \mu^{1/3}b_{e} .
     \end{align*}
     This is a conformal symplectic transformation with a multiplier $\mu^{-2/3}$, 
     then let $\mu\rightarrow 0$, consider the limiting case,
     the Hamiltonian for the Hill's lunar problem with the oblateness of the second primary is
     obtained as
     \begin{align}
      H_{b} &=\frac{1}{2}\|\mathbf{y}\|^{2}-(x_{1}y_{2}-x_{2}y_{1})-\frac{1}{r}
          -P_{2}(\frac{x_{1}}{r})r^{2}-\tilde{U} , \nonumber   \\ 
      \tilde{U} &= \frac{1}{r}\sum_{n=1}^{\infty}\left(\frac{a_{e}}{r}\right)^{2n}C_{2n} P_{2n}(\frac{x_{3}}{r}) .    
     \end{align}
     The Hill's lunar problem is studied by both analytical and numerical methods, see
     Michalodimitrakis (1980)\cite{Michalodi1980}, 
     Howision \& Meyer (2000I, 2000II)\cite{Howison00I,Howison00II},
     Maciejewski \& Rybicki (2001)\cite{Maciejewski2001}, 
     Llibre \& Roberto (2011)\cite{Llibre2011}, Belbruno et.al.(2019)\cite{Belbruno2019}.
     Special attention is paid to the effects of the oblateness on the motion of the 
     infinitesimal body by Sharma(1990)\cite{Sharma1990}, 
     Vashkov'yak \& Teslenko (2000,2001) \cite{Vashkovyak2000,Vashkovyak2001},
     Markellos et.al.(2001)\cite{Markellos2001}, 
     Perdiou (2008)\cite{Perdiou2008}, Bustos et.al.(2018)\cite{Bustos2018}. 
     
     In the paper of X.B. Xu (2019)\cite{XingboXu}, 
     a family of doubly symmetric periodic orbits of lunar type in the spatial CRTBP is
     shown to exist. Following a similar way of proof, such a family still exist in the Hill's lunar problem
     with the second primary oblate. The paper is organized as follows. 
     In Sect. 2, orbital elements and canonical elements are introduced into the 
     scaled Hamiltonian system with a small parameter.
     In Sect. 3, the short periodic effects of the first-order perturbation terms are eliminated  
     by the Lie transform method.
     In Sect. 4, the doubly-symmetric periodic solution is introduced and the continuation is 
     given by the use of a fixed point theorem. 
     In the last section, several discussions are proposed.
     
     %     \begin{align}
%    & \mathscr{R}_{1}: (x_{1},x_{2},x_{3},y_{1},y_{2},y_{3},t)
%    \rightarrow  (x_{1},-x_{2},-x_{3},-y_{1},y_{2},y_{3},-t), 
%    \nonumber    \\
%     & \mathscr{R}_{2}: (x_{1},x_{2},x_{3},y_{1},y_{2},y_{3},t)
%    \rightarrow  (x_{1},-x_{2},x_{3},-y_{1},y_{2},-y_{3},-t), 
%   \end{align}

     \section{The Hamiltonian in mixed elliptical elements}
        
       Set $\mathbf{x}\rightarrow \varepsilon^{2}\xi$,  
         $\mathbf{y}\rightarrow \varepsilon^{-1}\eta$, $r \rightarrow \varepsilon^{2}r$, 
         $a_{e} \rightarrow \varepsilon^{2}a_{e}$, 
         $C_{2n}\rightarrow \varepsilon^{6n}\tilde{J}_{2n}$, 
         this is a symplectic transformation with a multiplier $\varepsilon^{-1}$, and 
         the small parameter $\varepsilon$ represents the closeness of the infinitesimal body
         to the second primary. The Hamiltonian with $\varepsilon$ is             
   \begin{align}
     H_{c} &=\varepsilon^{-3}\left(\frac{\|\eta\|^{2}}{2}-\frac{1}{\|\xi\|}\right)-(\xi_{1}\eta_{2}-\xi_{2}\eta_{1})
     -\varepsilon^{3}r^{2}\mathcal{P}_{2}(\frac{\xi_{1}}{r})
     \nonumber  \\
    &  -\varepsilon^{-3}\frac{1}{\|\xi\|}\sum_{n=1}^{\infty}\left(\frac{a_{e}}{\|\xi\|}\right)^{2n}
     \varepsilon^{6n}\tilde{J}_{2n}P_{2n}(\frac{\xi_{3}}{\|\xi\|})  .
   \end{align}
   The perturbation caused by the oblateness is supposed to be as small as the 
   third-body perturbation from infinity, so the Hamiltonian $H_{c}$ 
   can be splitted into four parts,
   \begin{align}
     -H_{c}=\varepsilon^{-3}F_{01}(\xi,\eta)+F_{02}(\xi,\eta)+\varepsilon^{3}F_{1}(\xi)
     +\varepsilon^{9}F_{R}(\xi,\varepsilon) .
   \end{align}
   If the time is also scaled by $t=\varepsilon^{3}\tau$, then the Hamiltonian is converted
   to $\varepsilon^{3}H_{c}$. As $\varepsilon$ is small, $-F_{01}$ is a Kepler problem,
   so the Hamiltonian $\varepsilon^{3}H_{c}$ is a perturbed Keplerian system.
   
      The relations between the rectangular coordinates $(\vec{r},\dot{\vec{r}})$
       and the instantaneous orbital elements of a Keplerian orbit can refer to 
       some fundamental books \cite{Klioner,Duriez} or some papers 
       \cite{XingboXu, Xingbo2009} in celestial mechanics.
     The orbital elements $a,e,i$, $\Omega,\omega,M$, $f$ 
     are the semiaxis, eccentricity, inclination,
    the longitude of the ascending node, the argument of the pericenter, 
    the mean anomay and the eccentric anomaly, respectively.  
   Then the Delaunay elements, 
   \begin{align}
     \begin{array}{ccc}
       \mathcal{L}=\sqrt{a}, & \mathcal{G}=\sqrt{a(1-e^{2})}, & \mathcal{H}=\mathcal{G}\cos i,    \\
       h=\Omega,               &  g=\omega,                              & \ell =M ,
       \end{array} 
   \end{align}
   and the Poincar\'{e}-Delaunay elements, 
 \begin{align}
   \begin{array}{ccccccccc}
         Q_{1} &=&  \ell+g+h, &  Q_{2}  &=&  -\sqrt{2(\mathcal{L}-\mathcal{G})}\sin(g+h), 
    &  Q_{3} &=&  \ell+g,  \nonumber  \\
    \mathcal{P}_{1} &=&  \mathcal{L}-\mathcal{G}+\mathcal{H},   
    &  \mathcal{P}_{2}  &=&  \sqrt{2(\mathcal{L}-\mathcal{G})}\cos(g+h),  
    &  \mathcal{P}_{3} &=&   \mathcal{G}-\mathcal{H}.
    \end{array}
   \end{align}
  Parts of the Hamiltonian $H_{d}$ can be written as
   \begin{align}
    F_{01} &=\frac{1}{2\mathcal{L}^{2}}=\frac{1}{2(\mathcal{P}_{1}+\mathcal{P}_{3})^{2}}, \quad
    F_{02}=\mathcal{H}=\mathcal{P}_{1}-\frac{\mathcal{P}_{2}^{2}+Q_{2}^{2}}{2},   \nonumber \\
    F_{1} &=r^{2} P_{2}(\frac{\xi_{1}}{r})+
    \tilde{J}_{2n}\frac{a_{e}^{2}}{r^{3}} P_{2}(\frac{\xi_{3}}{r}),     
   \end{align}
   but $F_{1}$ cannot be expressed in a finite form of these canonical elements 
   because it contains the eccentric anomaly or the true anomaly, which is a function of the eccentricity 
   $e(\mathcal{L},\mathcal{G})$ and the mean anomaly $M=\ell$. 
   The first-order perturbation terms in $F_{1}$ can be expanded 
   by Tisserand expansion \cite{LaskarBoue}, one has 
   \begin{align}
     P_{2}(\frac{\xi_{1}}{r})
    &= I_{1}+I_{2}\cos 2(f+\omega-\Omega)+ I_{3}\cos 2(f+\omega+\Omega)
    \nonumber  \\
    & +I_{4}[\cos 2(f+\omega)+\cos 2\Omega] ,    
   \end{align} 
   where 
     \begin{align}
        I_{1}=\frac{3}{8}\frac{\mathcal{H}^{2}}{\mathcal{G}^{2}}-\frac{1}{8}, 
        I_{2}=\frac{3}{16}(1-\frac{\mathcal{H}}{\mathcal{G}})^{2},
        I_{3}=\frac{3}{16}(1+\frac{\mathcal{H}}{\mathcal{G}})^{2}, 
        I_{4}=\frac{3}{8}(1-\frac{\mathcal{H}^{2}}{\mathcal{G}^{2}}) ,
   \end{align}
   and one gets
   \begin{align}
      P_{2}(\frac{\xi_{3}}{r}) = -\frac{3}{4}\cos 2(f+\omega)+I_{5},
   \quad I_{5}=\frac{1}{4}-\frac{3}{4}\frac{\mathcal{H}^{2}}{\mathcal{G}^{2}}  .
    \end{align}
    According to the formulas of the Hansen coefficients \cite{LaskarBoue},
    one has 
    \begin{align*}
        \left(\frac{r}{a}\right)^{n}\exp (\texttt{i} m f)
        =\sum_{k=-\infty}^{\infty}X_{k}^{n,m}(e)\exp (\texttt{i} kM)    ,     
    \end{align*}
    so $r^{2}P_{2}(\frac{\xi_{1}}{r})$ can be expanded as
      \begin{align}\label{Hans1}
       r^{2}P_{2}(\frac{\xi_{1}}{r}) 
      &=I_{1}\sum_{k=-\infty}^{\infty}X_{k}^{2,0}(e)\cos(k Q_{1}-k(g+h))     \nonumber   \\
      &+  I_{2}\sum_{k=-\infty}^{\infty}X_{k}^{2,2}(e)\cos((k-4)Q_{1}+4Q_{3}+(2-k)(g+h))
       \nonumber  \\
      &  +I_{3}\sum_{k=-\infty}^{\infty}X_{k}^{2,2}(e)\cos(kQ_{1}+(2-k)(g+h))   \nonumber    \\
      & +I_{4}\sum_{k=-\infty}^{\infty}X_{k}^{2,2}(e)\cos((k-2)Q_{1}+2Q_{3}+(2-k)(g+h)) 
      \nonumber   \\
      &  +I_{4}\sum_{k=-\infty}^{\infty}X_{k}^{2,0}(e)\cos(kQ_{1}-k(g+h))\cos 2(Q_{1}-Q_{3}) , 
      \end{align}
      and $r^{-3} P_{2}(\frac{\xi_{3}}{r})$ can be expanded as
      \begin{align}
      r^{-3} P_{2}(\frac{\xi_{3}}{r})&=
      -\frac{3}{4}\sum_{k=-\infty}^{\infty}X_{k}^{-3,2}(e)
      \cos((k-2)Q_{1}+2Q_{3}+(2-k)(g+h))     \nonumber     \\
      & +I_{5}\sum_{k=-\infty}^{\infty}X_{k}^{-3,0}(e)\cos(kQ_{1}-k(g+h)) .
      \end{align}
      In addition, $g+h$ can be gotten from $Q_{2}$ and $\mathcal{P}_{2}$.
      The $e$ and $\cos i$ can be expressed as
      \begin{align}\label{eecosi}
       e^{2} =1-\left(1-\frac{\mathcal{P}_{2}^{2}+Q_{2}^{2}}
       {2(\mathcal{P}_{1}+\mathcal{P}_{3})}\right)^{2},   \quad 
      \frac{\mathcal{G}}{\mathcal{H}}=1-\frac{\mathcal{P}_{3}}{\mathcal{P}_{1}+\mathcal{P}_{3}
      -\frac{\mathcal{P}_{2}^{2}+Q_{2}^{2}}{2}} .
      \end{align}
      The Hansen coefficients above are very small if $|k|$ is large and $e$ small, 
       so only finite terms are needed to match the precision. 
      As there are infinite short-period terms, the first-order perturbation system 
      is complicate and can be simplified by the averaging method.
        
      \section{Averaging in the first-order system}
      
      Averaging method can be used to eliminate short periodic terms in 
      the perturbed dynamical system. In celestial mechanics, 
      Lie transforms is an explicit near-identy canonical transformation,
      and is taken in use here. There are two variables $Q_{1}, Q_{3}$
      which change as fast as the time because they both contain the mean anomaly. 
      Two averaging procedures are 
      needed in order to average $F_{1}$ over $Q_{1}$ and $Q_{3}$ in a period 
      of $2\pi$ successively.
      
      One has
     \begin{align}\label{aveF1}
    \bar{F}_{1} &=\frac{1}{2\pi}\int_{0}^{2\pi} F_{1} d Q_{1}
    \nonumber   \\
     & =  \biggl [ I_{1}X_{0}^{2,0}+I_{2}X_{4}^{2,2}\cos(4Q_{3}-2(g+h))+I_{3}X_{0}^{2,2}\cos(2g+2h)
     \nonumber  \\
     & +I_{4}X_{2}^{2,2}\cos 2Q_{3}+\frac{1}{2}I_{4}(X_{-2}^{2,0}+X_{2}^{2,0})\cos(2Q_{3}-2g-2h)   \biggr]
     \nonumber  \\
     & +\tilde{J}_{2} a_{e}^{2}\biggl[-\frac{3}{4}X_{2}^{-3,2}\cos 2Q_{3} +I_{5}X_{0}^{-3,0} \biggr],
   \end{align}
   and
   \begin{align}
     \bar{\bar{F}}_{1} & =\frac{1}{2\pi}\int_{0}^{2\pi}\bar{F}_{1}d Q_{3}  
     \nonumber  \\
    & = \left[ (1+\frac{3}{2}e^{2}) I_{1}+\frac{5}{2}  e^{2} I_{3} \cos(2g+2h)\right]
      +\tilde{J}_{2}a_{e}^{2}(1-e^{2})^{-3/2}I_{5} .
   \end{align}
   It is seen that $\bar{\bar{F}}_{1}$ does not contain terms about $\ell$ any more, 
   so $\mathcal{L}$ would be a constant in the first-order system, 
   and this makes it easier to give error estimates for 
   the continuation obeying to a fixed point theorem.
   
   Set $\tilde{\epsilon}=\varepsilon^{3}$ and let $\tilde{\epsilon}$ be the small parameter in the Lie transforms.
   The doubly averaged Hamiltonian is
   \begin{align}\label{Hdav}
       -\tilde{\epsilon}H_{d}=F_{01}+\tilde{\epsilon} F_{02} 
       +\frac{\tilde{\epsilon}^{2}}{2}\cdot 2\bar{\bar{F}}_{1}+O(\tilde{\epsilon}^{4}) , 
   \end{align}
   As there are two averaging procedures,  there are two generating functions 
   \begin{align}
      &   W=\frac{\tilde{\epsilon}^{2}}{2}\cdot 2W_{2} , \quad W_{2}=W_{2}^{(1)}+W_{2}^{(2)}, 
          \nonumber  \\
       &   \bar{F}_{1}-F_{1}+ \bar{\bar{F}}_{1}-\bar{F}_{1}
       =-\frac{\partial F_{01}}{\partial \mathcal{P}_{1}}\frac{\partial W_{2}^{(1)}}{\partial Q_{1}} 
         -\frac{\partial F_{01}}{\partial \mathcal{P}_{3}}\frac{\partial W_{2}^{(2)}}{\partial Q_{3}} ,    
   \end{align}
   the $W_{2}^{(1)}$ and $W_{2}^{(2)}$ can be achived by integrations,
   \begin{align}
    W_{2}^{(1)} &=(\mathcal{P}_{1}+\mathcal{P}_{3})^{3}\int \left(\bar{F}_{1}- F_{1}\right) dQ_{1},
         \nonumber  \\
    W_{2}^{(2)} &=(\mathcal{P}_{1}+\mathcal{P}_{3})^{3}\int \left(\bar{\bar{F}}_{1}- \bar{F}_{1}\right) dQ_{3},    
   \end{align}
   Neglecting the constants, one has
   \begin{align}
     -(\mathcal{P}_{1}+\mathcal{P}_{3})^{-3}W_{2}^{(1)} &=
     I_{1}\sum_{k \neq 0}\frac{1}{k}X_{k}^{2,0}\sin(kQ_{1}-k(g+h))    \nonumber   \\
     &+ I_{2}\sum_{k\neq 4}\frac{1}{k-4}X_{k}^{2,2}
     \sin((k-4)Q_{1}+4Q_{3}+(2-k)(g+h))                                             \nonumber   \\
     & +I_{3}\sum_{k\neq 0}\frac{1}{k}X_{k}^{2,2}\sin(kQ_{1}+(2-k)(g+h))   \nonumber  \\
     & +I_{4}\sum_{k\neq 2}\frac{1}{k-2}X_{k}^{2,2}\sin((k-2)Q_{1}+2Q_{3}+(2-k)(g+h))   \nonumber  \\
     & +\frac{1}{2}I_{4}\sum_{k\neq -2}\frac{1}{k+2}X_{k}^{2,0}\sin((k+2)Q_{1}-k(g+h)-2Q_{3}) 
     \nonumber  \\
     & +\frac{1}{2}I_{4}\sum_{k\neq 2}\frac{1}{k-2}X_{k}^{2,0}
     \sin((k-2)Q_{1}-k(g+h)+2Q_{3}), 
   \end{align}
    \begin{align}
     -(\mathcal{P}_{1}+\mathcal{P}_{3})^{-3}W_{2}^{(2)} &=
     \frac{1}{4}I_{2}X_{4}^{2,2}\sin(4Q_{3}-2(g+h))+\frac{1}{2}I_{4}X_{2}^{2,2}\sin 2Q_{3}  
     \nonumber  \\
     & +\frac{1}{4}I_{4}\left(X_{-2}^{2,0}+X_{2}^{2,0}\right)\sin(2Q_{3}-2(g+h))
     -\frac{3}{2}\tilde{J}_{2}a_{e}^{2}X_{2}^{-3,2}\sin 2Q_{3} .
    \end{align}
    Two generating functions are given above, and there are two times of Lie transforms.
    Both generating functions can be truncated according to a high order of the eccentricity.
    Theoretically, Lie transforms are inversible and can be calculated with the help of symbol
    calculation software and numerical calculation software \cite{Deprit1969, RandArmb}. 
   
   \section{Continuation of the doubly-symmetric periodic orbits}
    The truncated Hamiltonian $H_{appr}=-\tilde{\epsilon}^{-1}F_{01}-F_{02}$ is integrable, and is set as the 
    approximated system of the full system (\ref{Hdav}). 
    The aim of this section is to give the outline of the proof on the continuation of the 
    doubly symmetric periodic solutions.    
    
      There exists a lemma about the proposition about the doubly symmetric periodic solution, 
  \begin{lemma}\label{LemmaDoublyP}
    Consider an one-order autonomous ordinary differential system in $\mathbb{R}^{6}$, and it 
    is invariant under two anti-symplectic reflections:
      \begin{align}
  \mathscr{R}_{1} &: (x_{1},x_{2},x_{3},y_{1},y_{2},y_{3},t)\rightarrow
  (x_{1},-x_{2},-x_{3},-y_{1},y_{2},y_{3},-t),   \nonumber  \\
  \mathscr{R}_{2} &: (x_{1},x_{2},x_{3},y_{1},y_{2},y_{3},t)\rightarrow
  (x_{1},-x_{2},x_{3},-y_{1},y_{2},-y_{3},-t) .
 \end{align}
  That is to say, the system is symmetric about two Lagrangian planes, 
   \begin{align} 
   \mathscr{L}_{1} &=\{Z |Z=(x_{1},0,0,0,y_{2},y_{3})^{T}\},  \nonumber \\
   \mathscr{L}_{2} &=\{Z |Z=(x_{1},0,x_{3},0,y_{2},0)^{T}\} .
   \end{align}   
   If one solution hits the two Lagrangian planes $\mathscr{L}_{1}$ and $\mathscr{L}_{2}$ successively 
   with a time interval $T>0$, 
   then this solution is periodic with period $4T$ and doubly symmetric. 
   \end{lemma}
     
     These two Lagrangian planes can be expressed in Poincar\'{e}-Delaunay elements
     $Z=(Q_{1}$, $Q_{2}$, $Q_{3}$, $\mathcal{P}_{1},\mathcal{P}_{2},\mathcal{P}_{3})^{T}$,
      \begin{align}\label{PoinSymmetry}
    \mathscr{L}_{a} &=\{Z=(i \pi,  0,  j \pi, 
    \mathcal{P}_{1},\mathcal{P}_{2},\mathcal{P}_{3})^{T}  \},      
        \nonumber \\
    \mathscr{L}_{b} &=\{Z=(i \pi+k \pi, 0, j\pi+m\pi+\frac{\pi}{2}, \mathcal{P}_{1},
    \mathcal{P}_{2}, \mathcal{P}_{3})^{T}  \},
     \end{align}

    Denote the initial solution of the approximated system 
    as $Z_{0}^{\ast}\in \mathscr{L}_{a}$, 
    the differential equations are
      \begin{align}
      \begin{array}{llll}
      \frac{dQ_{1}}{d t}= \frac{\varepsilon^{-3}}{(\mathcal{P}_{1}+\mathcal{P}_{3})^{3}}-1, &
      \frac{dQ_{2}}{d t}= \mathcal{P}_{2}, &
      \frac{dQ_{3}}{d t}= \frac{\varepsilon^{-3}}{(\mathcal{P}_{1}+\mathcal{P}_{3})^{3}},     \\
      \frac{d\mathcal{P}_{1}}{d t}=0, &
      \frac{d\mathcal{P}_{2}}{d t}=- Q_{2}, &
      \frac{d\mathcal{P}_{3}}{d t}=0 .
      \end{array}
         \end{align}
      If $Z(t;Z_{0}^{\ast})$ is a doubly-symmetric periodic solution of the integrable approximated system,
      one has 
      \begin{align}
            Z_{0}^{\ast}=(i \pi,  0,  j \pi, \mathcal{P}_{1}^{\ast}, 0, \mathcal{P}_{3}^{\ast})^{T},
      \end{align}
      and one fourth of the period is $T_{0}^{\ast}=(m+\frac{1}{2}-k)\pi$. 
      
     Consider that the initial solution of the full system (\ref{Hdav}) belongs to $\mathscr{L}_{a}$, 
     and is near $Z_{0}^{\ast}$, that is
      \begin{align}\label{Z0f}
      Z_{0} &=Z_{0}^{\ast}+(0, 0, 0, \delta_{P1}, \delta_{P2}, \delta_{P3})^{T}
           \nonumber    \\
            &=(i\pi, 0,  j\pi, \mathcal{P}_{1}^{\ast}+\delta_{P1}, 0+\delta_{P2}, 
            \mathcal{P}_{3}^{\ast}+\delta_{P3})^{T},
      \end{align}    
      where $\delta_{P1}$, $\delta_{P2}$, $\delta_{P3}$ are small.
      
      The full system corresponds to the following kind of differential system,  
      \begin{align}\label{dotZ}
          \frac{dZ}{dt}=\mathcal{F}(Z,\tilde{\epsilon})=\tilde{\epsilon}^{-1}\mathcal{F}_{01}(Z)
          +\mathcal{F}_{02}(Z)
          +\tilde{\epsilon} \mathcal{F}_{1}(Z)
          +\tilde{\epsilon}^{2 +\kappa}\mathcal{F}_{R}(Z,\tilde{\epsilon}) ,
      \end{align}  
      where $Z\in \mathbb{R}^{n}$, $\kappa \geq 0$. The differential system
      of the approximated system
       \begin{align}
          \frac{dZ}{dt}=\mathcal{F}_{0}(Z,\tilde{\epsilon})=\tilde{\epsilon}^{-1}\mathcal{F}_{01}(Z)
                  +\mathcal{F}_{02}(Z)
        \end{align}
        is integrable, and has an analytical solution 
      $Z^{(0)}(t;z_{0},\tilde{\epsilon})$. 
      
      Consider that $\tilde{\epsilon}$ is sufficiently small, such that
      the difference between the solution of the full system and that of the approximated system 
      $\|Z(t;\mathcal{Z}_{0},\tilde{\epsilon})-Z^{(0)}(t;\mathcal{Z}_{0},\tilde{\epsilon})\|$ 
      remains small enough for a finite time interval. 
      The norm for vectors represents the maximum absolute value of its components.
       \begin{lemma}\label{lem2}
        Let $Z^{(1)}(t;z_{0},\tilde{\epsilon})$ be a solution of 
       \begin{align}\label{ZoneDiff}
            \dot{Z}^{(1)}(t;z_{0},\tilde{\epsilon})=\left(\mathbf{D}_{Z}\mathcal{F}_{0}\right)|_{Z^{(0)}} Z^{(1)}
            +\mathcal{F}_{1}(Z^{(0)})
       \end{align}
   with an initial condition $Z^{(1)}(0;z_{0},\tilde{\epsilon})=\mathbf{0}$.
    In a finite time interval $t\in [0,t_{0}]$, the solution of the full system (\ref{dotZ}) can be expressed as 
     \begin{align}\label{ZZ0Z1ZR}
      Z(t;z_{0},\tilde{\epsilon})=Z^{(0)}(t;z_{0},\tilde{\epsilon})
           +\tilde{\epsilon} Z^{(1)}(t;z_{0},\tilde{\epsilon})+\tilde{\epsilon}^{2}Z_{R}(t;z_{0},\tilde{\epsilon}),
     \end{align}
     where $\|Z(t;z_{0},\epsilon)-Z^{(0)}(t;z_{0},\tilde{\epsilon})\|$ is of order $\tilde{\epsilon}$.    
      The maximum absolute values for elements in
      $(\mathbf{D}_{Z}Z^{(1)})|_{Z_{0}} $ and $(\mathbf{D}_{Z}Z_{R})|_{Z_{0}}$
       are of zeroth order of $\tilde{\epsilon}$.       
      \end{lemma}

      Following Lemma \ref{lem2}, one can get the formulas of $Q_{1}$, $Q_{2}$,
      $Q_{3}$ in $\mathscr{L}_{b}$ with the initial solution $Z_{0}$ in (\ref{Z0f})
      for the full system after a finite time $T=T_{0}^{\ast}+\delta T =$ $(m+\frac{1}{2}-k)\pi $ $+\delta T$,
      \begin{align}\label{QQQ321}
    &  Q_{1}(T)=\left[\frac{\tilde{\epsilon}^{-1}}{(\mathcal{L}^{\ast}+\delta \mathcal{L})^{3}}-1 \right] T +i\pi
    +\tilde{\epsilon} Q_{1}^{(1)}+O(\tilde{\epsilon}^{2})=k\pi+i\pi,      \nonumber   \\
    &  Q_{2}(T)=\delta_{P2}\sin (T_{0}^{\ast}+\delta T) +\tilde{\epsilon} Q_{2}^{(1)}
            +O(\tilde{\epsilon}^{2})=0,     \nonumber   \\
    &  Q_{3}(T)=\frac{\tilde{\epsilon}^{-1}}{(\mathcal{L}^{\ast}+\delta \mathcal{L})^{3}}T+j\pi
    +\tilde{\epsilon} Q_{3}^{(1)}+O(\tilde{\epsilon}^{2})=j\pi+(m+\frac{1}{2})\pi ,
   \end{align}
     the equations can be reduced, firstly let the third equation minus the first equation, 
     secondly substitute 
     \begin{align}
     \tilde{\epsilon}^{-1}=\varepsilon^{-3}
     =\frac{(m+\frac{1}{2})}{m+\frac{1}{2}-k} (\mathcal{L}^{\ast})^{3}
     = \frac{(m+\frac{1}{2})\pi }{T_{0}^{\ast}} (\mathcal{L}^{\ast})^{3} , 
     \end{align} 
     into the third equation,
     finally one has
     \begin{align}
     \left\{\begin{array}{ll}
    &   \Psi_{1}= \delta T +\tilde{\epsilon} (Q_{3}^{(1)}-Q_{1}^{(1)})+O(\tilde{\epsilon}^{2})=0,                   
       \\
   &   \Psi_{2}=  (-1)^{m-k} \delta_{P2} \cos \delta T 
    +\tilde{\epsilon} Q_{2}^{(1)}+O(\tilde{\epsilon}^{2})=0,      \\
   &   \Psi_{3}=\frac{1+\frac{\delta T}{T_{0}^{\ast}}}
   {\left(1+\frac{\delta \mathcal{L}}{\mathcal{L}^{\ast}}\right)^{3}}-1
      +\frac{1}{(m+\frac{1}{2})\pi}\left[ \tilde{\epsilon} Q_{3}^{(1)}+\mathcal{O}(\tilde{\epsilon}^{2})\right]=0 .  
     \end{array}\right.
     \end{align}
     The three equations combine a vector $\Psi(\mathbf{X}, \tilde{\epsilon})$, where 
      $\mathbf{X}=(\delta T, \delta_{P2}, \delta \mathcal{L})^{T}$.  The Jacobian matrix
     derived from the partial derivatives of the above three equations over $\mathbf{X}$ is 
     non-degenerated, and one has
     \begin{align}
         \frac{\partial \Psi(\mathbf{X},0)}{\partial \mathbf{X}}=
          \left(\begin{array}{ccc}
             1                                                             &   0                   &    0     \\
     (-1)^{m-k+1}\sin \delta T \cdot \delta_{P2}        &   (-1)^{m-k}\cos \delta T     &    0     \\
         \frac{1}{T_{0}^{\ast}}\frac{1}{\left(1+\frac{\delta \mathcal{L}}{\mathcal{L}^{\ast}}\right)^{3} }
            &  0                   &   -\frac{3\left(1+\frac{\delta T}{T_{0}^{\ast}}\right)}
            {\left(1+\frac{\delta \mathcal{L}}{\mathcal{L}^{\ast}}\right)^{4}\mathcal{L}^{\ast}} 
        \end{array}\right) .
     \end{align}
     
     There is a corollary of Arenstorf's theorem given by Cors et.al. \cite{Cors05},
     \begin{lemma}[Cors, Pinyol \& Soler]\label{CorsSProp1}
       Let $\mathbb{U}$ be an open domain in $\mathbb{R}^{n}$, 
       $\mathbb{I}\subset \mathbb{R}$ an open neighbourhood of the 
    origin and $\mathbf{f}: \mathbb{U}\times \mathbb{I} \rightarrow \mathbb{R}^{n}$ 
    with $\mathbf{f}(\mathbf{0},0)=\mathbf{0}$, differentiable with respect to
    $x\in \mathbb{U}$, and $\mathbf{f}_{x}(\mathbf{0},0)$ non-singular. 
    Assume that there exist $c_{1}>0, c_{2}>0$ such that for $x\in \mathbb{U}$,
    $\epsilon \in \mathbb{I}$, 
    \begin{enumerate}
     \item    $\|\mathbf{f}_{x}(x,\epsilon)-\mathbf{f}_{x}(\mathbf{0},0)\|\leq  
      c_{1}(\|x\|+\epsilon)$, 
    \item    $\|\mathbf{f}(\mathbf{0},\epsilon)\|\leq  c_{2}\epsilon$ .
    \end{enumerate}
    Then there exists a function $x(\epsilon) \in \mathbb{U}$, 
    defined for $\epsilon \in \mathbb{I}' \subset \mathbb{I}$,
    such that $\mathbf{f}(x,\epsilon)=0$ and $x(0)=0$. 
     \end{lemma}
     Suppoe $\tilde{\epsilon}$ is in a neighborhood of zero and positive,
     one has 
     \begin{align}
       & \| \Psi_{\mathbf{X}}(\mathbf{X},\tilde{\epsilon})- \Psi_{\mathbf{X}}(\mathbf{0},0)\| \nonumber  \\
     & \leq \|\Psi_{\mathbf{X}}(\mathbf{X},\tilde{\epsilon})-\Psi_{\mathbf{X}}(\mathbf{X},0)\|+
     \| \Psi_{\mathbf{X}}(\mathbf{X},0)-\Psi_{\mathbf{X}}(\mathbf{0},0)\|
     \nonumber  \\
     &  \leq C_{1} \tilde{\epsilon}+ 
     C_{2}\max(\|\delta T\|,\|\delta_{P2}\|,\|\delta \mathcal{L}\|)
     \leq C_{3}(\|X\|+\tilde{\epsilon}),      
     \end{align}     
     and also
     \begin{align} 
          0<\|\Psi (\mathbf{0},\tilde{\epsilon}) \|<C_{4} \tilde{\epsilon} ,
      \end{align}
      where $C_{j}$ ($j=1,2,3,4$) are constants greater than zero. 
      One can choose the values of $\delta_{P1}$ and $\delta_{P2}$,
      under the condition of $\delta \mathcal{L}=\delta_{P1}+\delta_{P2}$. 
      After the continuation of $\mathbf{X}$, 
      the averaged initial value $Z_{0}$ can be transformed back to the original full system.
      The conclusion of this paper is
      \begin{theorem}
       For the spatial Hill's lunar problem with the second primary oblate,
       there exists a class of doubly-symmetric and near-circular
       periodic solutions around the oblate 
       primary. These orbits are symmetric with respect to the line joining two primaries,
      and to a plane. This plane contains that line connecting two primaries,
      and this plane is perpendicular to the primaries' motion plane. 
      \end{theorem}
      
      \section{Discussion}
      
      This paper completes a proof on the existence of a class of 
      doubly-symmetric and spatial near-circular 
      periodic solutions in the Hill's lunar problem with the second primary oblate.
      The method is almost the same as X.B.Xu(2019)\cite{XingboXu}. 
      New questions are arised after finishing this paper.
      For exampler, the stability and the global bifurcations of these orbits are still unknown,
       one can apply the averaging transforms, the shooting method and the 
      Poincar\'{e} cross section method to give some inspirations.  
      It will be interesting to calculate these solutions with the background of 
      astronomy in the future.

          \subsection*{Acknowledgements}
     The author would like to thank the reviewer of this paper for the comment.
      This work is supported by the National Nature Science Foundation 
     of China (NSFC, Grant No. 11703006).   
       
%\newpage
%and more a new one
%\newpage

\bigskip

\begin{flushleft}

{\bf AMS Subject Classification: ?????, ?????}\\[2ex]

% Write more than one author separately if they have different 
% affiliations, otherwise write the names on the same line, separeted 
% by commas.
%
Name~Xu,\\
Faculty of Mathematics and Physics, Institution of Huaiyin Technology \\
No. 1 of Meicheng Road, 223002 HuaiAn City, China\\
e-mail: \texttt{xbxu@hyit.edu.cn}\\[2ex]

% leave it blanck, you dont know these infos yet
%
\textit{Lavoro pervenuto in redazione il MM.GG.AAAA.}

\end{flushleft}

\normalsize
\label{\thechapter:lastpage}

\end{document}